# Exponentiable locales, revisited

Xu Huang

*Abstract.* We give a moderately motivated exposition of exponentiable locales and the construction of exponentials in Loc, without assuming prior knowledge of exponential topological spaces or continuous posets.

## 1  Introduction

Given two spaces $X$ and $Y$, a natural question is how to construct a space $X^Y$ representing the continuous functions from $Y$ to $X$. Fixing $Y$, if such a space can be constructed for each $X$, then we say $Y$ is *exponentiable*. The exponentiable objects in various categories of spaces, such as topological spaces and locales, have been thoroughly investigated. Escardó and Heckmann [1] gave an elementary account for the case of general topological spaces, and Hyland [2] discussed the case of locales.

The theory of exponentiable topological spaces proceeds by using the exponential law $\hom(1, X^Y) \cong \hom(Y, X)$, which determines that the points of this space are the continuous functions. We can then consider the coarsest or finest topology on this set making relevant functions continuous. The results can then be naturally modified to obtain localic versions. It is however difficult to see how one would arrive at such a result without the aid of topological spaces, since the universal property of exponentials naturally favors the characterization of their points.

This paper attempts to give a purely localic — or pointless — account of exponentials, without requiring prior knowledge of exponentiable topological spaces. Such an approach can then be applied to toposes, as detailed in Johnstone [3]. The goal is a sufficiently well-motivated exposition, plausible as a route by which exponential locales might have been discovered.

Ideally, this construction would result from pure "calculation", i.e. equational reasoning of natural isomorphisms. It is unclear whether this is achievable. Though if background knowledge of directed-complete partial orders is assumed, then there is an elegant and category-theoretic desciption by Townsend [4].

## 2  Preliminaries

We briefly state some definitions and fix notations.

### 2.1 Locales

**Definition 1.** A **frame** is a poset $F$ with finite meets and arbitrary joins, satisfying distributivity
$$a \wedge \bigvee_{\alpha \in I} b_\alpha = \bigvee_{\alpha \in I} a \wedge b_\alpha.$$

In particular, it contains the nullary meet $\top$ and nullary join $\bot$. A frame homomorphism is a map preserving finite meets and arbitrary joins.

**Definition 2.** A **locale** $X$ is defined by a frame $\Omega(X)$ whose elements are called **opens** of the locale. A **continuous map** $f : X \to Y$ is given by a frame homomorphism $f^* : \Omega(Y) \to \Omega(X)$. We denote its right adjoint — which exists by preservation of all joins — as $f_* : \Omega(X) \to \Omega(Y)$.

**Definition 3.** The **one-point space** 1 is defined by $\Omega(1) = \{\bot < \top\}$, or constructively the frame of propositions. A **point** of a locale $X$ is given by a continuous map $1 \to X$.

Intuitively, an open represents a yes-no question about the position of a point such that when the answer is yes, it is feasible to recognize this. For example, if the length of an object is within $(0.9, 1.1)$, then with sufficiently precise measurements, we can eventually contain the error within this interval. It would be infeasible to recognize when the length is *outside* the interval, for if it sits at $0.9$ precisely, no amount of measurement can tell us for sure. A point is then *defined* by a set of (consistent) answers to every such question.

**Definition 4.** The **Sierpinski space** $\mathbb{S}$ is defined by $\Omega(\mathbb{S}) = \{\bot < \omega < \top\}$, with $\omega$ denoting the *generic open*. Constructively, $\Omega(\mathbb{S})$ is the collection of upwards closed subsets of $\{0, 1\}$, with $\omega = \{1\}$. Continuous maps $X \to \mathbb{S}$ are in bijective correspondence with opens in $X$.

The most convenient way to present the product of two locales $X \times Y$, or the coproduct of two frames $\Omega(X) \otimes \Omega(Y)$, is via generators and relations. The frame is generated by tensor products $x \otimes y$, formally representing a rectangle, under the relations
$$(x \wedge x') \otimes (y \wedge y') = (x \otimes y) \wedge (x' \otimes y'),$$
$$\left(\bigvee_{\alpha \in I} x_\alpha\right) \otimes \left(\bigvee_{\beta \in J} y_\beta\right) = \bigvee_{\alpha \in I} \bigvee_{\beta \in J} x_\alpha \otimes y_\beta.$$
This is very similar to the notion of tensor products in algebra. [5]

However, in practice it is more convenient to characterize the tensor product using the machinery of coverages and $\mathcal{C}$-ideals. More specifically, we consider the set-theoretic product $\Omega(X) \times \Omega(Y)$, equipped with componentwise finite meets, and suggestively write $x \otimes y$ for the pair $(x, y)$. An element of $\Omega(X) \otimes \Omega(Y)$ corresponds to a downwards-closed subset $\varphi$ of $\Omega(X) \times \Omega(Y)$, such that if $\varphi$ contains $x_\alpha \otimes y_\beta$ for each $\alpha$ and $\beta$, then it contains $(\bigvee x_\alpha) \otimes (\bigvee y_\beta)$. Such a $\varphi$ is called a $\mathcal{C}$-ideal, and it formally represents the union of all the rectangles it contains, being *saturated* in the sense that all the geometrically compatible rectangles are already added to the union.

Given a downwards-closed subset $\psi$, we can characterize the $\mathcal{C}$-ideal $\langle \psi \rangle$ it generates. Naturally, we can consider

$$\psi' = \Big\{ x \otimes y \,\Big|\, x = \bigvee x_\alpha, y = \bigvee y_\beta, \quad \forall \alpha, \beta, x_\alpha \otimes y_\beta \in \psi \Big\},$$

i.e. closing $\psi$ under the coverage conditions once. Notably, $\psi' = \psi'' = \cdots = \langle \psi \rangle$. The fact that performing the construction once already suffices will be extremely useful in our analysis. This is the decategorified version of the plus construction in sheafification, which satisfies $P^{++} = P^{+++} = \cdots = a(P)$ for every presheaf $P$, stabilizing at the second step.

## 2.2 Exponentiable objects

**Definition 5.** Given a category $\mathcal{C}$ and two objects $X, Y$, suppose there is an object $E$ with a natural isomorphism
$$\hom(-, E) \cong \hom(- \times X, Y),$$
then we say $E$ (equipped with the natural isomorphism) is the **exponential object** $Y^X$. The image of the identity $\mathrm{id} : E \to E$ under this isomorphism is called the **evaluation map** $\mathrm{ev} : Y^X \times X \to Y$.

If $Y^X$ exists for all $Y$, then we say $X$ is **exponentiable**, in which case $(-)^X$ arranges into a functor, in particular the right adjoint to the functor $(- \times X)$.

# 3 Pursuing Exponentials

Our task is to figure out a criterion for when the exponential object in the category of locales $\mathsf{Loc}$ exists, and construct such an object.

The first reduction we can perform comes from the fact that $\Omega(\mathbb{S})$ is the free frame generated by one generator $\omega$. As is the case in all algebraic structures, every frame can be presented by generators and relations, which can be expressed as a colimit of copies of $\Omega(\mathbb{S})$. Dualizing, every locale is a limit of copies of $\mathbb{S}$. Since the right adjoint $(-)^X$ preserves limits, we can first focus our attention to $\mathbb{S}^X$.

Using the universal property of exponentials, we can easily see that the points of $\mathbb{S}^X$ are exactly the opens of $X$, since $\hom(1, \mathbb{S}^X) \cong \hom(X, \mathbb{S})$, but it doesn't really help us because locales are far from determined by its points, a priori. To obtain the set of opens, we can try to compute $\hom(\mathbb{S}^X, \mathbb{S}) \cong \Omega(\mathbb{S}^X)$, but there are no obvious bijections allowing us to do that. However, we can at least get some obvious opens in this way: given a point $p : 1 \to X$ of $X$, we have a map $\mathbb{S}^X \to \mathbb{S}^1 \cong \mathbb{S}$, which corresponds to an open $O_p$.

Another lead we have is the adjoint functor theorem. Naïvely, the theorem states that if $F$ is a functor preserving colimits, then $F$ has a right adjoint $G$ given by the formula
$$G(X) = \operatorname*{colim}_{F(Y) \to X} Y. \tag{$\star$}$$
Therefore, we can potentially obtain a criterion of exponentiability by studying when $(- \times X)$ preserves colimits, and a construction of the exponential by simplifying the limit formula.

The coproduct of locales is dually the product of frames. Since frames are algebraic structures, the product is given by the set-theoretic product equipped with componentwise operations. It is not difficult to see that tensor products distribute over arbitrary direct products in Frm, which proves that $(- \times X)$ always preserves arbitrary coproducts in Loc. What remains is to characterize when the functor preserves coequalizers, or equivalently when tensor products preserve equalizers in Frm.

The situation has some semblance the study of tensor products of modules, where tensor products always distribute over direct sums, and a module is *flat* when it preserves equalizers.

On the other hand, using $(\star)$ we can obtain an expression of $\mathbb{S}^X$:
$$\mathbb{S}^X = \underset{u \in \Omega(Y \times X)}{\text{colim}} Y$$
where the colimit is over the category of pairs $(Y, u)$, with arrows $(Y, u) \to (Z, w)$ given by continuous maps $f : Y \to Z$ such that $u = f^*(w)$. Note that the partial order of the opens is irrelevant to the limit diagram. Dualizing,
$$\Omega(\mathbb{S}^X) = \underset{u \in \Omega(Y) \otimes \Omega(X)}{\lim} \Omega(Y).$$
Again, since frames are algebraic, the limit can be calculated by a limit in Set, equipped with componentwise operations. Ignoring size issues, an element $\alpha$ of the limit set is given by a family of maps $\alpha_Y : \Omega(Y) \otimes \Omega(X) \to \Omega(Y)$ (without preserving any structure) *natural in $Y$*. Intuitively, we have an element of $\Omega(Y) \otimes \Omega(X)$, and we want to consider ways to extract an open in $Y$ that depends "uniformly" in $Y$.

As an example, given any point $p : 1 \to X$, we can take $\hat{p} : Y \times 1 \to Y \times X$ and consider $\hat{p}^*(u) \in \Omega(Y \times 1) \cong \Omega(Y)$. This gives a natural family of maps, which corresponds to an element of $\Omega(\mathbb{S}^X)$. Indeed, this is the open $O_p \in \Omega(\mathbb{S}^X)$ mentioned earlier.

**Remark.** It is also possible to obtain this characterization via the Yoneda lemma. An open of $\mathbb{S}^X$ is given by a map $\mathbb{S}^X \to \mathbb{S}$, which is determined by a natural transformation $\hom(Y, \mathbb{S}^X) \to \hom(Y, \mathbb{S})$ in $Y$, which is exactly a natural family of maps $\Omega(Y \times X) \to \Omega(Y)$.

## 3.1 Natural operations

Since points are not fundamental in the study of locales, we should consider an open $s \in \Omega(X)$ and see if there is a corresponding natural map $\Omega(Y) \otimes \Omega(X) \to \Omega(Y)$ associated with it. Geometrically, we can draw the following picture:

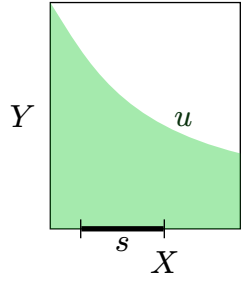

In the previous example, we have a point $p : 1 \to X$, and our map $\hat{p}^*$ simply slices the open $u \in \Omega(Y) \otimes \Omega(X)$ vertically at $p$ and takes the intersection. For the case of an open $s \in \Omega(X)$, the natural thought is to consider the largest rectangle $y \otimes s$ contained in $u$, and take $y$ as the output. Algebraically, this corresponds to expressing $u$ as a sum of tensors

$$u = \bigvee_{s \in \Omega(X)} y_s \otimes s$$

and taking the coefficient $y_s$. (Since there may be multiple ways to express $u$ as a sum, we take the one with largest coefficient to make it well-defined.) In other words, we are considering an operation $\text{coeff}_s : \Omega(Y) \otimes \Omega(X) \to \Omega(Y)$ for each $s$ given by

$$\text{coeff}_s(u) = \bigvee \{y \in \Omega(Y) \,|\, y \otimes s \leq u\}.$$

and we want to verify the naturality of this construction, i.e. for any continuous map $f : Z \to Y$, we need

$$f^* \text{coeff}_s(u) \stackrel{?}{=} \text{coeff}_s\left(\hat{f}^* u\right)$$

where $\hat{f} : Z \times X \to Y \times X$ is the product map.[1]

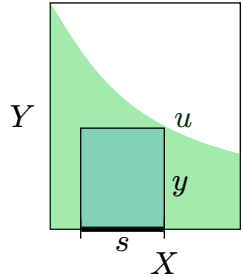

It's easy to see the left hand side is contained in the right hand side: since $\text{coeff}_s(\hat{f}^* u)$ is the largest coefficient $z \in \Omega(Z)$ such that the rectangle $z \otimes s$ is contained in $\hat{f}^*(u)$, and whenever $y \otimes s \leq u$ we have $f^*(y) \otimes s \leq \hat{f}^*(u)$, by definition we have $f^*(y) \leq \text{coeff}_s(\hat{f}^* u)$ for every such $y$, and therefore $f^* \text{coeff}_s(u) \leq \text{coeff}_s(\hat{f}^* u)$.

For the other direction, we need to use the $\mathcal{C}$-ideal characterization of $\Omega(Y) \otimes \Omega(X)$. In particular, we need to show

$$\text{coeff}_s\left(\hat{f}^* u\right) \leq f^* \text{coeff}_s(u)$$
$$\iff \bigvee \left\{z \in \Omega(Z) \,\big|\, z \otimes s \leq \hat{f}^*(u)\right\} \leq f^* \text{coeff}_s(u)$$
$$\iff \left[\forall z, z \otimes s \leq \hat{f}^*(u) \implies z \leq f^* \text{coeff}_s(u)\right].$$

---

[1]This is also why the other obvious thing — taking the meet of all the rectangles that contains the slice — doesn't work: infinite meets aren't preserved, even though they always exist.

How do we simplify $\hat{f}^*$? By definition it acts as $\hat{f}^*(y \otimes s) = f^*(y) \otimes s$ on the rectangles, so we need to find a way to express $u$ as a sum of rectangles. Since we have

$$u = \bigvee_{s \in \Omega(X)} \operatorname{coeff}_s(u) \otimes s \qquad (\dagger)$$

we can now deduce

$$\hat{f}^*(u) = \bigvee_{s \in \Omega(X)} f^* \operatorname{coeff}_s(u) \otimes s.$$

This is a union of rectangles, and we know it contains $z \otimes s$. This intuitively seems to imply $z \leq f^* \operatorname{coeff}_s(u)$, but it is not a priori true since the containment of $z \otimes s$ may have resulted from several other summands.

Using the $\mathcal{C}$-ideal characterization, $\hat{f}^*(u)$ is regarded as an ideal generated by $f^* \operatorname{coeff}_s(u) \otimes s$. By our discussion in Section 2.1, a rectangle $z \otimes s$ is contained in $\hat{f}^*(u)$ iff there exists some cover $z = \bigvee_\alpha z_\alpha$ and $s = \bigvee_\beta s_\beta$ such that each $z_\alpha \otimes s_\beta$ is contained in some rectangle $f^* \operatorname{coeff}_{s_{\alpha,\beta}}(u) \otimes s_{\alpha,\beta}$. Since $\operatorname{coeff}_s(u)$ is antitone in $s$, we can safely take $s_{\alpha,\beta} = s_\beta$, as this gives the best chance of making $z_\alpha \leq f^* \operatorname{coeff}_{s_{\alpha,\beta}}(u)$ hold.

To summarize the current proof state, we have $z = \bigvee_\alpha z_\alpha$ and $s = \bigvee_\beta s_\beta$ such that $z_\alpha \leq f^* \operatorname{coeff}_{s_\beta}(u)$, and we wish to prove $z \leq f^* \operatorname{coeff}_s(u)$. Without loss of generality, we can assume there is only one $\alpha$, because in the general case, we can apply the special case for each $\alpha$ to get $\forall \alpha, z_\alpha \leq f^* \operatorname{coeff}_s(u)$ which implies $z \leq f^* \operatorname{coeff}_s(u)$ anyway. Hence the goal is to prove that

$$z \leq f^* \operatorname{coeff}_{s_\beta}(u) \implies z \leq f^* \operatorname{coeff}_s(u)$$

whenever $\bigvee_\beta s_\beta = s$.

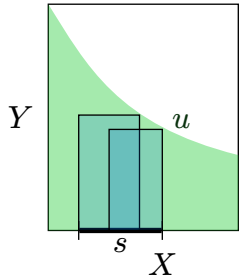

This is equivalent to

$$\bigwedge_\beta \operatorname{coeff}_{s_\beta}(u) \leq \operatorname{coeff}_s(u).$$

which seems plausible from the picture, but the infinite meet should raise some suspicion. Substituting the definition of $\operatorname{coeff}_s(u)$, we have

$$\bigwedge_\beta \bigvee \{y \in \Omega(Y) \mid y \otimes s_\beta \leq u\} \leq \bigvee \{y \in \Omega(Y) \mid y \otimes s \leq u\}. \qquad (\ddagger)$$

Suppose there were finitely many $\beta$, then I claim this inequality would be true. We can distribute the meet into the join on the left, which gives us $\bigvee_{y_\beta \otimes s_\beta \leq u} \bigwedge_\beta y_\beta$ where the join ranges over all families $y_\beta$. So we just need to prove $\bigwedge_\beta y_\beta \leq \bigvee \{y \mid y \otimes s \leq u\}$ for each such $y_\beta$. Again, since there are finitely many $\beta$, we have an identity

$$\left(\bigwedge_\beta y_\beta\right) \otimes \left(\bigvee_\beta s_\beta\right) \leq \bigvee_\beta (y_\beta \otimes s_\beta). \tag{¶}$$

But since the right hand side is contained in $u$ by assumption and $\bigvee_\beta s_\beta = s$, we see that $(\bigwedge_\beta y_\beta) \otimes s \leq u$ and $\bigvee_\beta y_\beta$ is a summand in the right hand side of (‡). Therefore the inequality holds *if the intersection is finite*.

Indeed, there are cases where naturality fails. Consider $X = [0, +\infty)$, $Y = [0, 1]$ is the interval and $Z = \{0\}$ is a one-point subspace of $Y$. Let $u$ be the open under the curve $y = \exp(-u)$.

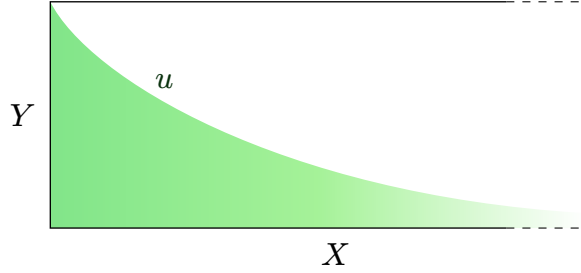

Suppose $s = \top \in \Omega(X)$ is the entire space. There is no open of $Y$ such that $y \otimes s \leq u$. However, if we let $f : Z \to Y$ be the subspace inclusion, then $\hat{f}^*(u)$ is the entire $x$-axis, so
$$\text{coeff}_s\left(\hat{f}^* u\right) = \top$$
$$f^* \text{coeff}_s(u) = \bot$$
which don't agree. If we follow along the proof, we will see that
$$\bigvee_{n=0}^\infty [0, e^{-n}) \otimes [0, n) \leq u$$
and indeed $[0, n)$ is a (genuinely) infinite cover of $[0, +\infty)$.

## 3.2 Proof repair

If we reflect on this counterexample, we see that the reason this fails can be explained by the intuition of opens. As mentioned in Section 2.1, an open should correspond to a question about the location of a point, such that when the answer is yes, it should be feasible to verify it. A point in $\mathbb{S}^X$ is an open of $X$, and so an open of $\mathbb{S}^X$ should be a question about open sets $o \in \Omega(X)$. Translating our $\text{coeff}_s$ construction via the adjoint functor theorem, this corresponds to the question

Is the open $s$ contained in $o$? [2]

But this is not feasible for verification: if $o$ is "touching the boundary" of $s$, then a tiny error would change the answer of the question. This happens when $o = (0, 1)$ and $s = \left(\frac{1}{2}, 1\right)$ in $\mathbb{R}$. Of course, topologically $(0, 1)$ is indistinguishable from $(0, +\infty)$, so the the boundary at infinity also counts, as shown in the counterexample.

---

[2] This can be verified by taking $Y = 1$ and following the proof of the naïve adjoint functor theorem.

With this in mind, we can consider an alternative construction. Suppose there is some to-be-determined relation $s \ll o$ formalizing the intuition of $s$ being well-contained in $o$ with room for error, then the question "Is the open $s$ well-contained in $o$" can be seen as replacing

$$\operatorname{coeff}_s(u) = \bigvee \{y \in \Omega(Y) \mid y \otimes s \leq u\}$$

with

$$F_s(u) = \bigvee \{y \in \Omega(Y) \mid \exists s' \gg s,\ y \otimes s' \leq u\}$$

We can replay the argument above with most $\operatorname{coeff}_s$ replaced by $F_s$. (Note that (†) still uses $\operatorname{coeff}_s$ instead of $F_s$.) Intuitively, we should require $s \ll o$ to imply $s \leq o$, and this guarantees the first half of the argument, i.e.

$$f^* F_s(u) \leq F_s\big(\hat{f}^* u\big)$$

still holds. The second half then reduces to proving

$$\bigwedge_\beta \operatorname{coeff}_{s_\beta}(u) \leq F_s(u).$$

Here we have $s \ll s' = \bigvee_\beta s_\beta$ replacing the original $s = \bigvee_\beta s_\beta$. Recall that if there were only finitely many $\beta$, then the proof would work. Suppose $\beta \in I$. We shall take a finite subset $J \subseteq I$ and investigate what kinds of subsets are required. We have

$$\bigwedge_{\beta \in I} \operatorname{coeff}_{s_\beta}(u) \leq \bigwedge_{\beta \in J} \operatorname{coeff}_{s_\beta}(u)$$

so it suffices to prove $\bigvee_{\beta \in J} \operatorname{coeff}_{s_\beta}(u) \leq F_s(u)$. Going through the proof and using the identity (¶), we see that the requirement for $J$ is that $s \ll \bigvee_{\beta \in J} s_\beta$, making the reasonable assumption that if $s' \leq s \ll o \leq o'$ then $s' \ll o'$ too.

Summarizing, we need a relation $s \ll o$ such that if there is a cover $\bigvee_{\beta \in I} o_\beta = o$, then we have a finite subset $J \subseteq I$ such that $\bigvee_{\beta \in J} o_\beta \gg s$ still holds. Presumably, not every locale $X$ can support such a relation, for not every locale is exponentiable. Hence, we will need to figure out both a definition of $\ll$ and a criterion of exponentiability, so that they jointly imply this condition.

### 3.3 Exponentiability

As explained at the beginning of this section, exponentiability of locales hinges on the existence of $\mathbb{S}^A$. Let us state and prove this rigorously.

**Theorem 6.** A locale $A$ is exponentiable if and only if the exponential $\mathbb{S}^A$ exists.

*Proof.* One direction is trivial. Suppose $\mathbb{S}^A$ exists and $Y$ is an arbitrary locale. A canonical presentation of the frame $\Omega(Y)$ is given by having a generator $\omega_y$ for each $y \in \Omega(Y)$, and set

$$\begin{aligned} \bigwedge_{\beta \in J} \omega_{y_\beta} &= \omega_{\bigwedge_{\beta \in J} y_\beta} \\ \bigvee_{\alpha \in I} \omega_{y_\alpha} &= \omega_{\bigvee_{\alpha \in I} y_\alpha} \end{aligned} \qquad (\divideontimes)$$

where $J$ is finite and $I$ is arbitrary, and $y_\alpha$, $y_\beta$ are families of elements without repeat. This means we the collection of all such equalities forms a set $K$. For $k \in K$ we write $\lambda_k$ for the left hand side and $\rho_k$ for the right hand side.

Take a coproduct $F_0$ of $\Omega(Y)$-many copies of $\Omega(\mathbb{S})$, and another coproduct $F_1$ of $K$-many copies of $\Omega(\mathbb{S})$. They are isomorphic to the free frames generated by the corresponding sets. Now $\lambda_k$ and $\rho_k$ each provide a map $K \to F_0$, which extends to a map $F_1 \to F_0$. The coequalizer of these maps are then canonically isomorphic to $\Omega(Y)$. Thus $\Omega(Y)$ is a nested colimit of $\Omega(\mathbb{S})$, and hence $Y$ as a nested limit of $\mathbb{S}$.

Now suppose we have a diagram $D : \mathcal{J} \to \mathsf{Loc}$ such that $D(i)^A$ exists for each $i \in \mathcal{J}$. Since $\mathsf{Loc}$ is complete, we can calculate
$$\hom\left(-, \lim_{i \in \mathcal{J}} D(i)^A\right)$$
$$\cong \lim_{i \in \mathcal{J}} \hom\left(-, D(i)^A\right)$$
$$\cong \lim_{i \in \mathcal{J}} \hom(- \times A, D(i))$$
$$\cong \lim_{i \in \mathcal{J}} \hom(- \times A, D(i))$$
$$\cong \hom\left(- \times A, \lim_{i \in \mathcal{J}} D(i)\right).$$
Therefore $\lim_i D(i)^A$ satisfies the universal property of the exponential $(\lim_i D(i))^A$. This shows that the full subcategory spanned by $Y$ such that $Y^A$ exists is closed under limits. Therefore if $\mathbb{S}$ is in such this subcategory, then every object is. □

**Remark.** Colimits of algebraic structures with known presentations are also easy to present with generators and relations. This makes Theorem 6 very effective.

Focusing our attention on $\mathbb{S}^A$, we would like to obtain a criterion on $A$ — more precisely a criterion on the structure of $\Omega(A)$ — such that $\mathbb{S}^A$ exists. Since the opens of $A$ correspond to the points of $\mathbb{S}^A$, this seems to be our only angle of attack.

To clarify the exact relationship between them, given an open $a \in \Omega(A)$, we have a locale map $\tilde{p}_a : A \to \mathbb{S}$ defined by $\tilde{p}_a^*$ mapping the generic open $\omega$ to $a$. This is then transposed to a map $p_a : 1 \to \mathbb{S}^A$. On the other hand, given such a point, we can recover $\tilde{p}_a$ by the composition
$$A \xrightarrow{\sim} 1 \times A \xrightarrow{p_a \times \mathrm{id}} \mathbb{S}^A \times A \xrightarrow{\mathrm{ev}} \mathbb{S}.$$

These are all formal properties that hold for all categories. How do we obtain something specific to $\mathsf{Loc}$? We don't yet have anything to say about $\mathbb{S}^A$, but we do know how products are constructed. We need a continuous map $\mathrm{ev} : \mathbb{S}^A \times A \to \mathbb{S}$, but this is completely specified by an open of $\mathbb{S}^A \times A$, which can be written as a union of rectangles. Formally,
$$\mathrm{ev}^*(\omega) = \bigvee_{\substack{d' \in \Omega(\mathbb{S}^A) \\ d' \otimes a' \leq \mathrm{ev}^*(\omega)}} \bigvee_{a' \in \Omega(A)} d' \otimes a'.$$

Therefore we can apply the product map $p_a \times \mathrm{id}$ and get
$$\top \otimes a = (p_a \times \mathrm{id})^* \mathrm{ev}^*(\omega)$$
$$= \bigvee_{\substack{d' \in \Omega(\mathbb{S}^A)\, a' \in \Omega(A) \\ d' \otimes a' \leq \mathrm{ev}^*(\omega)}} \bigvee (p_a \times \mathrm{id})^*(d' \otimes a')$$
$$= \bigvee_{\substack{d' \in \Omega(\mathbb{S}^A)\, a' \in \Omega(A) \\ d' \otimes a' \leq \mathrm{ev}^*(\omega)}} \bigvee p_a^*(d') \otimes a'.$$

The isomorphism between $A$ and $1 \times A$ is given by mapping $\varphi \otimes a$ to $\bigvee\{a \mid \varphi = \top\}$. Classically this takes $\varphi \otimes a$ to $a$ if $\varphi = \top$ and $\bot$ otherwise. In addition, $p_a^*(d') = \top$ is by definition $p_a \in d'$, where $p_a$ is a point and $d'$ an open. Using this we can further simplify
$$a = \bigvee_{\substack{d' \in \Omega(\mathbb{S}^A)\, a' \in \Omega(A) \\ d' \otimes a' \leq \mathrm{ev}^*(\omega)}} \bigvee \bigvee \{a' \mid p_a^*(d')\}$$
$$= \bigvee \{a' \in \Omega(A) \mid \exists d' \in \Omega(\mathbb{S}^A), d' \otimes a' \leq \mathrm{ev}^*(\omega) \text{ and } p_a \in d'\}.$$

On first glance, this seems too complicated for any use. However, there are still useful information to be extracted. In particular, since the summands of a join are by definition contained in the join, this expression reveals a small fact:

**Lemma 7.** *If there is an open $d' \in \Omega(\mathbb{S}^A)$ such that $p_a \in d'$ and $d' \otimes a' \leq \mathrm{ev}^*(\omega)$, then $a' \leq a$.*

Here, we are using the partial order $a' \leq a$ on the opens of $A$. Since we are also investigating $p_a : 1 \to \mathbb{S}^A$, we need to determine how the partial order translates to the points. In fact, there is a natural partial order on points, or more generally continuous maps, by comparing the frame homomorphism pointwise. Now we have three kinds of objects: opens $a \in \Omega(A)$, maps $\tilde{p}_a : A \to \mathbb{S}$ and points $p_a : 1 \to \mathbb{S}^A$. It's obvious by definition that the order on $\tilde{p}_a$ agrees with that on $a$. So the question is whether the exponential adjunction
$$\hom(Z, Y^X) \cong \hom(Z \times X, Y)$$
respects the partial order. The uncurrying map sends $F : Z \to Y^X$ to
$$Z \times X \xrightarrow{F \times \mathrm{id}} Y^X \times X \xrightarrow{\mathrm{ev}} Y.$$
But it's easy to verify that both function composition and product maps respect the order, so uncurrying is monotonic, which implies monotonicity for currying too.

Although this might seem mundane, the coincidence of these partial orders have consequences. Since elementhood of points $p \in u$ is defined as $p^*(u) = \top$, if $p \leq q$ in the partial order, we have $p^*(u) = \top \implies q^*(u) = \top$, and hence $p \in u \implies q \in u$. This is topologically known as the **specialization order** on points, since a point is more special if it falls in fewer open sets. Returning to our case, we see that open sets in $\mathbb{S}^A$ are upwards-closed, in the sense that if $d$ contains $p_a$, then it will contain $p_{a'}$ for $a \leq a'$.

It should be warned that although the order on $A \to \mathbb{S}$ and $1 \to \mathbb{S}^A$ is defined pointwise, the joins and meets are not necessarily pointwise. We can see this is true for $A \to \mathbb{S}$, but there is no reason for it to be true for $1 \to \mathbb{S}^A$. In other words,
$$p^*_{a \vee a'}(d) \neq p^*_a(d) \vee p^*_{a'}(d).$$
And we write $p_a \vee p_{a'}$ for the "true" join, i.e. the least upper bound in $1 \to \mathbb{S}^A$, which is equal to $p_{a \vee a'}$.

Viewing them as frame homomorphisms, a point $1 \to \mathbb{S}^A$ is given by the subset of $\Omega(\mathbb{S}^A)$ mapped to $\top$, and the pointwise joins and meets corresponds to the unions and intersections of this subset, which are no longer frame homomorphisms. This is similar to how the joins (i.e. least upper bounds) of subalgebras are not given by unions, because algebraic operations force new elements to belong to the join.

**Remark.** For subalgebras, meets are always computed as intersections. But in the case of points, the frame homomorphism requirement not only forces new elements to be in the join, but also forces some old elements to be evicted from the meet. So neither joins nor meets are pointwise.

### 3.4 Local compactness

Now we seem to be stuck on both fronts. On one hand, during the construction of $\mathbb{S}^A$ we would like a binary relation $a' \ll a$ such that if $a$ is covered by opens, then there is a subset whose union $a''$ still satisfies $a' \ll a''$. On the other hand, if $\mathbb{S}^A$ exists, then the opens of $A$ can be expressed as a strange union with complex conditions.

We have to make both ends meet somehow. So the only move is to consider how they interact. Recall that if $\mathbb{S}^A$ exists, then
$$a = \bigvee \{a' \in \Omega(A) \mid \exists d' \in \Omega(\mathbb{S}^A), d' \otimes a' \leq \mathrm{ev}^*(\omega) \text{ and } p_a \in d'\}. \tag{§}$$
and in particular $a' \leq a$ for all such $a'$. We can attempt to investigate covers $a = \bigvee_\alpha a_\alpha$ and the joins of its finite subsets. However, we just warned that the joins $\bigvee_\alpha p_{a_\alpha}$ cannot be computed pointwise, so it would be difficult to understand when $\bigvee_\alpha p_{a_\alpha} \in d'$ may hold. Nonetheless, there is a partial result that can be salvaged.

Inspired by the case of computing joins in subalgebras, we see that *directed* joins can still be computed pointwise. In other words, if given a family of points, for any two elements $p, q$, there is a third $r$ greater than both of them (and the family is non-empty), then — again using the correspondence of homomorphisms $\Omega(X) \to \Omega(1)$ with special subsets of $\Omega(X)$ for intuition — for every new element forced to join the subset via $p$ and $q$, such an element has already joined the subset via $r$.

Luckily, since we also care about the finite joins in the family, there is a natural way to turn the join into a directed join. Given a set $M$ of elements, we can consider the set $M^\uparrow$ of all finite joins of these points. This set is now directed, with the same join as $M$.

In our situation, we have a family of opens $a_\alpha \in \Omega(A)$, $\alpha \in I$, and we can consider the family of finite joins $a_J = \bigvee_{\alpha \in J} a_\alpha$ for each finite $J \subseteq I$. Now the join $\bigvee_J p_{a_J}$ of the points can be computed pointwise, and it is still equal to $p_a$. We have $p_a \in d'$, or equivalently $p_a^*(d') = \top$, so there must be a particular point $p_b$ in $M^\uparrow$ such that $p_b^*(d') = \top$, equivalently $p_b \in d'$. By definition, $b = a_J$ for some finite $J$.

Going back to our starting point, we have $a' \leq a$ whenever there is $d' \in \Omega(\mathbb{S}^A)$ such that $p_a \in d'$ and $d' \otimes a' \leq \mathrm{ev}^*(\omega)$. However, we just proved that $p_{a_J} \in d'$ as well, so we must conclude that $a' \leq a_J$ too. Summarizing, we have proved a strengthening of Lemma 7.

**Lemma 8.** Suppose we have $d' \in \Omega(\mathbb{S}^A)$ such that $p_a \in d'$ and $d' \otimes a' \leq \mathrm{ev}^*(\omega)$. Let $a = \bigvee_{\alpha \in I} a_\alpha$ be a cover, then there is a finite subset $J \subseteq I$ with $a' \leq \bigvee_{\alpha \in J} a_\alpha$.

The definition has now revealed herself.

**Definition 9.** Given two opens $s$ and $o$, $s \ll o$ ($s$ is **way below** $o$) iff for every covering of $o$, there is a finite subset of opens that covers $s$.

A nice consequence of this definition is that $s \ll s$ iff $s$ is compact. So this definition is a kind of *relative* compactness. It is indeed almost equivalent to the topological definition of relative compactness, but we shall leave this to the reader.

We can restate Lemma 8 more succintly as follows: if there is $d' \in \Omega(\mathbb{S}^A)$ such that $p_a \in d'$ and $d' \otimes a' \leq \mathrm{ev}^*(\omega)$, then $a' \ll a$. Notice that there is still some potence left in (§) that we have not extracted: namely that the join of all these $a'$ is exactly $a$. In particular, this means $a = \bigvee_{a' \ll a} a'$ for all $a \in \Omega(A)$ whenever $\mathbb{S}^A$ exists. We can produce from this a definition.

**Definition 10.** A locale is **locally compact** iff every open is the union of subopens way below it.

We can now verify the desired property of $\ll$.

**Lemma 11.** Given a locally compact locale and opens $a' \ll a$, if there is an open cover of $a$, then there is a finite subset whose join is way above $a'$.

*Proof.* We just need to find an open $s$ such that $a' \ll s \ll a$. Consider the set $\{t \mid \exists s, t \ll s \ll a\}$. Using local compactness twice we see this set covers $a$, and there is a finite subset $t_k \ll s_k \ll a$ such that $\bigvee_k t_k = a'$. Take $s = \bigvee_k s_k$, which satisfies $s \ll a$ since it is a finite join of elements way below $a$. Similarly, $a' \ll s$ since $a'$ is a finite join of elements way below $s$. This finishes the proof. $\square$

## 4 Main Result

Our pursuit has come to bear fruit:

**Theorem 12** *(Hyland).* The following are equivalent:
   (i) $A$ is exponentiable;

(ii) $\mathbb{S}^A$ exists;

(iii) $A$ is locally compact.

*Proof.* (i) $\implies$ (ii) is obvious. Our discussion above has shown that (ii) $\implies$ (iii). Earlier, we have also proved (ii) $\implies$ (i) as Theorem 6. What is left is to prove (iii) $\implies$ (ii). We will also give a concrete construction of the general exponential space.

By Lemma 11 we conclude that each map $F_s$ we constructed does indeed correspond to an open of $\mathbb{S}^A$, when $A$ is locally compact. As discussed, an open of $\mathbb{S}^A$ is a question about the position of a point $p_a$ in $\mathbb{S}^A$, equivalently an open $a$ in $A$. The open for $F_s$ corresponds to the question "Is $s$ well contained in $a$?" We can thus suggestively notate it as $\lceil s \ll \mathrm{O} \rceil$, where O denotes the putative open in question. This is to be regarded as a family of formal symbols indexed by $s$, which we claim form a set of generators for $\Omega(\mathbb{S}^A)$.

To finish the construction of $\Omega(\mathbb{S}^A)$, we just to pose some relations that $\lceil s \ll \mathrm{O} \rceil$ should satisfy. These relations are either obvious, or will suggest themselves when we attempt to prove $\mathbb{S}^A$ is indeed an exponential. However, we list them here for convenience. These can either be phrased as logical implications of the questions $\lceil s \ll \mathrm{O} \rceil$, or as the containment of opens.

- If $s \leq s'$, then $\lceil s' \ll \mathrm{O} \rceil$ implies $\lceil s \ll \mathrm{O} \rceil$, i.e. the former open contains the latter.
- For a *finite* join $\bigvee_\alpha s_\alpha = s$, the conjunction of $\lceil s_\alpha \ll \mathrm{O} \rceil$ implies $\lceil s \ll \mathrm{O} \rceil$.
- $\lceil s \ll \mathrm{O} \rceil$ implies the disjunction of $\lceil s' \ll \mathrm{O} \rceil$ where $s' \gg s$.

To show that $\mathbb{S}^A$ is an exponential, we need to construct a natural isomorphism
$$\hom(-, \mathbb{S}^A) \cong \hom(- \times A, \mathbb{S}).$$
From right to left, we have an open $u \in \Omega(Z \times A)$, and we wish to construct a continuous map $\hat{u} : Z \to \mathbb{S}^A$, which amounts to a map $\hat{u}^*$ sending $\lceil s \ll \mathrm{O} \rceil$ to $\Omega(Z)$ preserving the relations. By the geometric intuition, we set
$$\hat{u}^* \lceil s \ll \mathrm{O} \rceil = \bigvee \{z \mid \exists s' \gg s, z \otimes s' \leq u\}.$$
From left to right, we have a continuous map $F$ sending $\lceil s \ll \mathrm{O} \rceil$ to $F^* \lceil s \ll \mathrm{O} \rceil \in \Omega(Z)$, such that the relations are preserved, and we need to construct an open $\bar{F} \in \Omega(Z \times A)$ corresponding to the map $Z \times A \to \mathbb{S}$. Now the best we can do is
$$\bar{F} = \bigvee_{s \in \Omega(A)} F^* \lceil s \ll \mathrm{O} \rceil \otimes s,$$
so we would need local compactness to guarantee this recovers the open. In particular, the evaluation map $\mathrm{ev} : \mathbb{S}^A \times A \to \mathbb{S}$ is given by the open
$$\overline{\mathrm{id}} = \bigvee_{s \in \Omega(A)} \lceil s \ll \mathrm{O} \rceil \otimes s.$$
The proof now proceeds by verifying the two constructions are mutual inverses, and that they are natural in $Z$. These can proceed in a similar way to Section 3.1 and Section 3.2, using the machinery of $\mathcal{C}$-ideals. We leave this to the reader. $\square$

## 4.1 Explicit construction of general exponentials

We also record the results of applying Theorem 6 to obtain a construction of $\Omega(B^A)$ for arbitrary $B$. We start with a coproduct of $\Omega(B)$-many copies of $\Omega(\mathbb{S}^A)$. This is generated by the opens $\lceil s \ll O_b \rceil$ for each $s \in \Omega(A)$ and $b \in \Omega(b)$, where $b$ serves as a label for the copies. For a fixed $b$, the opens $\lceil s \ll O_b \rceil$ will satisfy the relations as posed in $\Omega(\mathbb{S}^A)$. We then need to work out the extra relations given by the frame coequalizer in Theorem 6, which will produce $\Omega(B^A)$.

Geometrically, an open of $B^A$ asks a question about the configuration of a function $f : A \to B$. After quotienting by the new relations, the image of $\lceil s \ll O_b \rceil$ in $\Omega(B^A)$ will correspond to the question "Is $s$ well contained in $f^*(b)$?" We can suggestively write this as $\lceil s \ll f^*(b) \rceil$ instead.

To calculate the relations on these generators, let's temporarily write the product of $\Omega(B)$-many copies of a locale $L$ as $L^{|\Omega(B)|}$.[3] Recall the relations came from (※), requiring equalities between pairs of opens $\Omega(\mathbb{S}^{|\Omega(B)|})$, or equivalently maps $\mathbb{S}^{|\Omega(B)|} \to \mathbb{S}$, to obtain $B$. After exponentiating, we have pairs of maps $(\mathbb{S}^A)^{|\Omega(B)|} \to \mathbb{S}^A$ that we want to equate, where we identified $(\mathbb{S}^A)^{|\Omega(B)|} \cong (\mathbb{S}^{|\Omega(B)|})^A$. For clarity, we write

$$\mathfrak{E} : \Omega(\mathbb{S}^{|\Omega(B)|}) \to \hom((\mathbb{S}^A)^{|\Omega(B)|}, \mathbb{S}^A)$$

for this conversion map.

However, it is not straightforward to determine the image of elements in (※) under $\mathfrak{E}$. If $\iota_b : L^{|\Omega(B)|} \to L$ are the projection maps, then plausibly $\mathfrak{E}(\omega_b) = \iota_b$, which we will see to be true. The joins and meets are much more difficult. We do not know if $\mathfrak{E}$ preserves the joins and meets, and even so, they cannot be computed pointwise in the poset $\hom((\mathbb{S}^A)^{|\Omega(B)|}, \mathbb{S}^A)$.

Hence, we need to clarify how $\mathfrak{E}$ is defined. $\mathfrak{E}(x)$ is given by taking the composition

$$(\mathbb{S}^A)^{|\Omega(B)|} \times A \xrightarrow{\sim} (\mathbb{S}^{|\Omega(B)|})^A \times A \xrightarrow{\text{ev}} \mathbb{S}^{|\Omega(B)|} \xrightarrow{\tilde{p}_x} \mathbb{S}$$

and then the exponential transpose. If $x = \omega_b$, then the corresponding map $\tilde{p}_x = \iota_b$, and the composition gives an open of $(\mathbb{S}^A)^{|\Omega(B)|} \times A$ defined by

$$u_b = \bigvee_{a \in \Omega(A)} \lceil a \ll O_b \rceil \otimes a.$$

The exponential transpose is given by

$$\mathfrak{E}(\omega_b)^* \lceil a \ll O \rceil = \bigvee \{ z \mid \exists a' \gg a, z \otimes a' \leq u_b \}, \tag{$*$}$$

where $z$ ranges over the opens of $(\mathbb{S}^A)^{|\Omega(B)|}$. We shall proceed to simplify this construction.

We can substitute and use the $\mathcal{C}$-ideal characterization of $u_b$ to get $z \otimes a' \leq u_b$ is equivalent to the existence of families $a'_\beta$ and $a''_\gamma$ such that $\lceil a'_\beta \ll O_b \rceil \leq \lceil a''_\gamma \ll O_b \rceil$ for all $\beta \in J$ and $\gamma \in K$, and

$$z \le \bigvee_{\beta \in J} \lceil a'_\beta \ll \mathrm{O}_b \rceil \quad \text{and} \quad a' \le \bigvee_{\gamma \in K} a''_\gamma.$$

---

[3]This is in fact homeomorphic to an exponentiation with discrete locales, but we will not need this.

Since $a' \gg a$, we can pick a finite subset $K' \subseteq K$ of $a''_\gamma$ that covers $a$. In this case we can simplify

$$\lceil a'_\beta \ll \mathrm{O}_b \rceil \le \bigwedge_{\gamma \in K'} \lceil a''_\gamma \ll \mathrm{O}_b \rceil \le \lceil a \ll \mathrm{O}_b \rceil$$

where the last inequality occurs in the presentation of $\Omega(\mathbb{S}^A)$. This means a necessary condition of $\exists a' \gg a, z \otimes a' \le u_b$ is $z \le \lceil a \ll \mathrm{O}_b \rceil$. This is almost a sufficient condition, but recall in Lemma 11 we showed that in a locally compact locale $a \ll a'$ implies the existence if an interpolating $a \ll a_0 \ll a'$. So in the discussion above, instead of picking $K' \subseteq K$ such that $a''_{\gamma \in K'}$ covers $a$, we can pick $a \ll a_0 \ll a'$ first, so that $a''_{\gamma \in K'}$ covers $a_0$ instead. This modified argument shows

$$\exists a' \gg a, z \otimes a' \le u_b \implies \exists a_0 \gg a, \ z \le \lceil a_0 \ll \mathrm{O}_b \rceil.$$

Using the definition of $u_b$ the reverse implication is also obvious.

Next, a basis of $L^{|\Omega(B)|}$ is given by finite intersections of the opens $\iota_b^*(l)$, with $b \in \Omega(B)$ and $l \in \Omega(L)$ — boxes in the product space with finitely many non-trivial factors. So for our purposes we may assume $z$ to be of the form $\bigwedge_{\alpha \in J} \lceil a_\alpha \ll \mathrm{O}_{b_\alpha} \rceil$ where $J$ is finite. We can then geometrically see that this intersection is contained in $\lceil a_0 \ll \mathrm{O}_b \rceil$ iff *either* one of the factors has $b_\alpha = b$ and $\lceil a_\alpha \ll \mathrm{O}_b \rceil \le \lceil a_0 \ll \mathrm{O}_b \rceil$ *or* one of the factors is empty. For the latter, $z$ is empty and doesn't contribute to the join in $(*)$. This means

$$\mathfrak{E}(\omega_b)^* \lceil a \ll \mathrm{O} \rceil = \bigvee \{ \lceil a' \ll \mathrm{O}_b \rceil \mid \exists a_0 \gg a, \lceil a' \ll \mathrm{O}_b \rceil \le \lceil a_0 \ll \mathrm{O}_b \rceil \}$$
$$= \bigvee \{ \lceil a_0 \ll \mathrm{O}_b \rceil \mid a_0 \gg a \}$$
$$= \lceil a \ll \mathrm{O}_b \rceil.$$

So we confirmed that $\mathfrak{E}(\omega_b) = \iota_b$ as expected. This is a lot of work to check something we intuitively know, but it lays the foundation for our next task, which is to calculate $\mathfrak{E}(x)$ for the other opens $x$ involed in $(\divideontimes)$.

Suppose $x \in \Omega(\mathbb{S}^{|\Omega(B)|})$ is given by a finite meet or arbitrary join of $\omega_b$'s, then the corresponding open of $(\mathbb{S}^A)^{|\Omega(B)|} \times A$ will be the meet or join of $u_b$'s, since frame homomorphisms preserve them. So, $(*)$ is still true if we replace $u_b$ by the corresponding meet or join.

For finite meets, it suffices to deal with the nullary and binary cases. In the nullary case,

$$\mathfrak{E}(\top)^* \lceil a \ll \mathrm{O} \rceil = \bigvee \{ z \mid \exists a' \gg a \} = \begin{cases} \top & \exists a' \gg a \\ \bot & \nexists a' \gg a \end{cases}$$

so in $\Omega(B^A)$ we have, for the first case,

$$\lceil a \ll f^*(\top) \rceil = \top \quad (a \ll \top)$$

since the existence of some $a' \gg a$ is equivalent to $\top \gg a$ by monotonicity. The second case ($\nexists a' \gg a$) is implied by the relation

$$\lceil a \ll f^*(b) \rceil = \bigvee_{a \ll a'} \lceil a' \ll f^*(b) \rceil$$

which we already have. (Constructively, imposing the first relation is also enough.)

For binary meets $u_{b_1} \wedge u_{b_2}$, we need to have a uniform $a' \gg a$ such that both $z \otimes a' \leq u_{b_1}$ and $z \otimes a' \leq u_{b_2}$ holds. Running through the argument, we have

$$\exists a_0 \gg a,\ z \leq \lceil a_1 \ll O_{b_1} \rceil \wedge \lceil a_2 \ll O_{b_2} \rceil$$

as a necessary condition for $\exists a' \gg a, z \otimes a' \leq u_{b_1} \wedge u_{b_2}$, which is obviously sufficient too. Again, $z$ can be taken as a finite meet of $\lceil a \ll O_b \rceil$ in $(*)$. A similar calculation shows

$$\mathfrak{E}(\omega_{b_1} \wedge \omega_{b_2})^* \lceil a \ll O \rceil$$
$$= \bigvee \left\{ \lceil a_1' \ll O_{b_1'} \rceil \wedge \lceil a_2' \ll O_{b_2'} \rceil \,\middle|\, \exists a_0 \gg a, \begin{matrix} \lceil a_1' \ll O_{b_1} \rceil \leq \lceil a_0 \ll O_{b_1} \rceil \\ \lceil a_2' \ll O_{b_2} \rceil \leq \lceil a_0 \ll O_{b_2} \rceil \end{matrix} \right\}$$
$$= \bigvee \left\{ \lceil a_0 \ll O_{b_1} \rceil \wedge \lceil a_0 \ll O_{b_2} \rceil \,\middle|\, a_0 \gg a \right\}.$$

This cannot be further simplified, so in $\Omega(B^A)$ our relation is

$$\lceil a_0 \ll f^*(b_1) \rceil \wedge \lceil a_0 \ll f^*(b_2) \rceil \leq \lceil a \ll f^*(b_1 \wedge b_2) \rceil.$$

Similar to the nullary case, we only need an inequality because the converse inequality is implied by another relation we already have, establishing the equality.

We are left with arbitrary joins. It is possible to directly proceed and calculate $\mathfrak{E}(\bigvee \omega_b)$, but we can reduce our workload by noting that joins of directed sets are, again, very easy to calculate. In fact, it's not hard to see

$$\mathfrak{E}\left( \bigvee_{\alpha \in I} \omega_{b_\alpha} \right) = \bigvee_{\alpha \in I} \iota_{b_\alpha}$$

and the join can be computed pointwise. Therefore our relation in $\Omega(B^A)$ is

$$\lceil a \ll f^*(b) \rceil = \bigvee_{\alpha \in I} \lceil a \ll f^*(b_\alpha) \rceil$$

for $\{b_\alpha\}$ a directed cover of $b$. Incidentally, this also implies $\lceil a \ll f^*(b) \rceil$ is monotonic in $b$, if we take the directed join to be over $\{b, b'\}$ where $b' \leq b$.

As discussed earlier, an arbitrary join can be decomposed into a directed join of finite joins. So our final goal is finite joins. Here, a similar argument shows

$$\mathfrak{E}\left( \bigvee_{\alpha \in I} \omega_{b_\alpha} \right)^* \lceil a \ll O \rceil$$
$$= \bigvee \left\{ z \,\middle|\, \exists a_\alpha \text{ covers } a, z \leq \bigwedge_{\alpha \in I} \lceil a_\alpha \ll O_{b_\alpha} \rceil \right\}$$
$$= \bigvee_{a_\alpha} \bigwedge_{\alpha \in I} \lceil a_\alpha \ll O_{b_\alpha} \rceil$$

and the corresponding relation is

$$\lceil a \ll f^*(b) \rceil = \bigvee_{a_\alpha} \bigwedge_{\alpha \in I} \lceil a_\alpha \ll f^*(b_\alpha) \rceil$$

where $a_\alpha$ ranges over finite covers of $a$.

We collect the complete presentation of $\Omega(B^A)$ in Figure 1.

$$\lceil a \ll f^*(b) \rceil \vdash \lceil a' \ll f^*(b') \rceil \qquad (a' \leq a, b \leq b')$$
$$\vdash \lceil \bot \ll f^*(b) \rceil$$
$$\lceil a \ll f^*(b) \rceil, \lceil a' \ll f^*(b) \rceil \vdash \lceil a \vee a' \ll f^*(b) \rceil$$
$$\vdash \lceil a \ll f^*(\top) \rceil \qquad (a \ll \top)$$
$$\lceil a \ll f^*(b) \rceil, \lceil a \ll f^*(b') \rceil \vdash \lceil a' \ll f^*(b \wedge b') \rceil \qquad (a' \ll a)$$
$$\lceil a \ll f^*(b) \rceil \vdash \bigvee_{a \ll a'} \lceil a' \ll f^*(b) \rceil$$
$$\lceil a \ll f^*(b) \rceil \vdash \bigvee_{\alpha \in I} \lceil a \ll f^*(b_\alpha) \rceil \qquad (\{b_\alpha\} \text{ directed cover of } b)$$
$$\lceil a \ll f^*(b) \rceil \vdash \bigvee_{a_\alpha} \bigwedge_{\alpha \in I} \lceil a_\alpha \ll f^*(b_\alpha) \rceil \qquad (\{b_\alpha\} \text{ finite cover of } b)$$

Figure 1. Hyland's axiomatization of $\Omega(B^A)$